\documentclass[10pt,thmsa]{article}
\usepackage{amsfonts}

\usepackage{sw20lart}
\usepackage[magyar]{babel}
\usepackage[T1]{fontenc}
\usepackage[latin2]{inputenc}

\input tcilatex
\newtheorem{Tetel}{Theorem}

\newtheorem{Lemma}[Tetel]{Lemma}

\newenvironment{bizonyitas}{{\bf Proof:}}{\hfill $\Box$}
\begin{document}

\title{Uniform convergence for convexification\\
of dominated pointwise convergent continuous functions}
\date{}
\author{Zoltán Kánnai \\
%EndAName
Department of Mathematics\\
Corvinus University Budapest}
\maketitle

\noindent \textbf{Abstract:} The Lebesgue dominated convergence theorem of
the measure theory implies that the Riemann integral of a bounded sequence
of continuous functions over the interval $\left[ 0,1\right] $ pointwise
converging to zero, also converges to zero. The validity of this result is
independent of measure theory, on the other hand, this result together with
only elementary functional analysis, can generate measure theory itself \cite
{daniell}. The mentioned result was also known before the appearance of
measure theory \cite{Arzela}\cite{Osgood}, but the original proof was very
complicated. For this reason this result, when presented in teaching, is
generally obtained based on measure theory. Later, Eberlein \cite{eberlein}
gave an elementary, but still relatively complicated proof, and there were
other simpler proofs but burdened with complicated concepts, like measure
theory (e.g. \cite{lewin}). In this paper we give a short and elementary
proof even for the following strenghened form of the mentioned result: a
bounded sequence of continuous functions defined on a compact topological
space $K,$ pointwise converging to zero, has a suitable convexification
converging also uniformly to zero on $K,$ thus, e.g., the original sequence
converges weakly to zero in $C\left( K\right) .$ This fact can also be used
in the proof of the Krein-Smulian theorem \cite{whitley}. The usual proof
beyond the simple tools of the functional analysis, uses heavy embedding
theorems and the Riesz' representation theorem with the whole apparatus of
measure theory. Our main result, however, reduces the cited proof to a form
in which we need abstract tools only, namely the Hahn-Banach separation
theorem and Alaoglu's theorem, without Riesz' representation or any
statement of measure theory.

\section{Positive functionals keep the dominated pointwise convergence}

Let $K$ be a compact topological space, denote by $C\left( K\right) $ the
Banach space of real valued continuous functions over $K.$ The operator
''upper envelope'' (resp. ''lower envelope'') is denoted by the symbol $\vee 
$ (resp. $\wedge $). Fix a positive linear functional $\Phi :C\left(
K\right) \rightarrow \Bbb{R}.$ For every function $f:K\rightarrow \Bbb{R}%
_{+} $, define 
\[
\overline{\Phi }\left( f\right) :=\sup\limits_{\QATOPD. . {g\in C\left(
K\right) }{0\leq g\leq f}}\Phi g\;. 
\]
$\overline{\Phi }$ is obviously monotone, i.e. $\overline{\Phi }\left(
f_{1}\right) \geq \overline{\Phi }\left( f_{2}\right) $ if $f_{1}\geq f_{2},$
moreover $\overline{\Phi }$ is an extension of $\Phi .$

We note that $\overline{\Phi }$ is not subadditive. This fact gives
importance to the next lemma.

\begin{lemma}
Let $f_{1},f_{2}:K\rightarrow \Bbb{R}_{+},$ $f_{1}\geq f_{2}$ and $g,h\in
C\left( K\right) ,$ $0\leq h\leq f_{1},$ $0\leq g\leq f_{2}.$ Then 
\[
\overline{\Phi }\left( f_{2}\right) -\Phi \left( g\wedge h\right) \leq
\left( \overline{\Phi }\left( f_{1}\right) -\Phi h\right) +\left( \overline{%
\Phi }\left( f_{2}\right) -\Phi g\right) \;.
\]
\end{lemma}

\begin{bizonyitas}
In view of the identity $g+h=\left( g\vee h\right) +\left( g\wedge h\right) $
we obtain 
\[
\Phi g+\Phi h=\Phi \left( g\vee h\right) +\Phi \left( g\wedge h\right) \;,
\]
hence $g\vee h\leq f_{1}$ implies 
\[
\Phi g+\Phi h\leq \overline{\Phi }\left( f_{1}\right) +\Phi \left( g\wedge
h\right) \;,
\]
whence obviously 
\[
\overline{\Phi }\left( f_{2}\right) -\Phi \left( g\wedge h\right) \leq
\left( \overline{\Phi }\left( f_{1}\right) -\Phi h\right) +\left( \overline{%
\Phi }\left( f_{2}\right) -\Phi g\right) \;.
\]
\end{bizonyitas}

\begin{lemma}
Take a sequence of functions $f_{n}:K\rightarrow \Bbb{R}_{+}$ with $%
f_{1}\leq \alpha <+\infty $ and $f_{n}\left( x\right) $ decreasingly tending
to $0$ for every $x\in K.$ Then $\overline{\Phi }\left( f_{n}\right)
\rightarrow 0.$
\end{lemma}

\begin{bizonyitas}
Fix a number $\varepsilon >0.$ For every integer $n$ choose a function $%
g_{n}\in C\left( K\right) $ such that $0\leq g_{n}\leq f_{n}$ and 
\[
\Phi g_{n}>\overline{\Phi }\left( f_{n}\right) -\frac{\varepsilon \;}{%
\;2^{n+1}}\;.
\]
Put $h_{n}:=\bigwedge\limits_{k=1}^{n}g_{k}.$ Prove by induction that 
\[
\overline{\Phi }\left( f_{n}\right) -\Phi h_{n}\leq \sum\limits_{k=1}^{n}%
\frac{\varepsilon \;}{\;2^{k+1}}\;.
\]
The case $n=1$ is obvious since $h_{1}=g_{1}.$ Supposing the statement is
true for $n,$ from the previous lemma we get 
\begin{eqnarray*}
\overline{\Phi }\left( f_{n+1}\right) -\Phi h_{n+1} &=&\overline{\Phi }%
\left( f_{n+1}\right) -\Phi \left( g_{n+1}\wedge h_{n}\right) \leq  \\
&\leq &\left( \overline{\Phi }\left( f_{n}\right) -\Phi h_{n}\right) +\left( 
\overline{\Phi }\left( f_{n+1}\right) -\Phi g_{n+1}\right) \leq  \\
&\leq &\sum\limits_{k=1}^{n}\frac{\varepsilon \;}{\;2^{k+1}}+\frac{%
\varepsilon \;}{\;2^{n+2}}=\sum\limits_{k=1}^{n+1}\frac{\varepsilon \;}{%
\;2^{k+1}}\;,
\end{eqnarray*}
which completes the induction. Now the sequence $\left( h_{n}\right)
\subseteq C\left( K\right) $ tends pointwise decreasingly to $0,$ so by
Dini's theorem we obtain that $h_{n}$ uniformly tends to $0$ on $K,$ hence
by continuity of $\Phi $ we get $\Phi h_{h}\rightarrow 0.$ Thus, there is an
index $N$ such that $0\leq \Phi h_{h}\leq \frac{\varepsilon }{2}$ for every $%
n\geq N,$ consequently $0\leq \overline{\Phi }\left( f_{n}\right) \leq \Phi
h_{n}+\sum\limits_{k=1}^{n}\frac{\varepsilon \;}{\;2^{k+1}}\leq \frac{%
\varepsilon }{2}+\sum\limits_{k=1}^{n}\frac{\varepsilon \;}{\;2^{k+1}}%
<\varepsilon .$
\end{bizonyitas}

\begin{theorem}
Let $\left( g_{n}\right) \subseteq C\left( K\right) $ with $0\leq g_{n}\leq
\alpha <+\infty $ for every $n\in \Bbb{N},$ pointwise converging to $0.$
Then $\Phi g_{n}\rightarrow 0.$
\end{theorem}

\begin{bizonyitas}
The sequence of functions $f_{n}:=\bigvee\limits_{k=n}^{\infty }g_{k}$
fulfills the assumptions of the previous lemma. Thus, 
\[
0\leq \Phi g_{n}\leq \overline{\Phi }\left( f_{n}\right) \rightarrow 0\;.
\]
\end{bizonyitas}

This theorem yields, for example, that a dominated sequence $\left(
g_{n}\right) \subseteq C\left[ 0,1\right] $ tending pointwise to $0$ has
also $\int\limits_{0}^{1}g_{n}\rightarrow 0,$ without any application of
measure theory, or other complicated ideas and methods.

\section{A convexification admitting uniform convergence}

Let $K$ be a compact topological space.

\begin{Lemma}
Let $\left( f_{n}\right) \subseteq C\left( K\right) $ with $0\leq f_{n}\leq
\alpha <+\infty $ for every $n\in \Bbb{N},$ pointwise converging to $0.$
Then for arbitrary $\varepsilon >0,$ there is a convex combination $g$ of
the functions $f_{1},f_{2},\ldots $ such that $g\leq \varepsilon $ on the
whole $K.$
\end{Lemma}

\begin{bizonyitas}
Suppose, by contradiction, that there is an $\varepsilon >0$ such that the
convex subsets 
\[
M:=\limfunc{co}\left( f_{1},f_{2},\ldots \right) \text{\textrm{\qquad
and\qquad }}\varepsilon +C\left( K\right) _{-}:=\left\{ f\in C\left(
K\right) :f\leq \varepsilon \right\} 
\]
of $C\left( K\right) $ are disjoint. Since the interior of $\varepsilon
+C\left( K\right) _{-}$ is nonempty, by the Hahn-Banach separation theorem
we obtain a nonzero functional $\Phi \in \left( C\left( K\right) \right) ^{*}
$ such that 
\[
\sup\limits_{\varepsilon +C\left( K\right) _{-}}\Phi \leq
\inf\limits_{M}\Phi \;.
\]
This implies that the functional $\Phi $ is positive. Then by the latest
theorem, $\Phi f_{n}\rightarrow 0,$ on the other hand 
\[
\Phi f_{n}\geq \inf\limits_{M}\Phi \geq \sup\limits_{\varepsilon +C\left(
K\right) _{-}}\Phi \geq \Phi \left( \varepsilon \cdot \mathbf{1}\right) >0
\]
for every $n,$ since $\Phi $ is positive and nonvoid. This is a
contradiction.
\end{bizonyitas}

\begin{definition}
Let $\left( f_{n}\right) $ be a sequence from a set $X.$ A sequnce $\left(
g_{n}\right) \subseteq X$ is called a \emph{convexification} of the sequence 
$\left( f_{n}\right) $ if 
\[
g_{n}\in \limfunc{co}\left( f_{n},f_{n+1},\ldots \right) 
\]
for every $n\in \Bbb{N}.$
\end{definition}

\begin{theorem}
\label{convexifi}Let $\left( f_{n}\right) \subseteq C\left( K\right) $ with $%
\left| f_{n}\right| \leq \alpha <+\infty $ for every $n\in \Bbb{N},$
pointwise converging to $0.$ Then there is a convexification $\left(
g_{n}\right) $ of the sequence $\left( f_{n}\right) $ converging uniformly
to $0$ on $K.$
\end{theorem}

\begin{bizonyitas}
Apply the previous lemma for $\varepsilon =\frac{1}{n}$ and the sequence $%
\left| f_{n}\right| ,\left| f_{n+1}\right| ,\ldots .$
\end{bizonyitas}

\begin{remark}
From the previous theorem immediately follows that an arbitrary functional $%
\Phi \in \left( C\left( K\right) \right) ^{*}$ has the property of dominated 
$0$-convergence on $C\left( K\right) ,$ without the Jordan-type
decomposition.
\end{remark}

\section{An application: Krein's theorem with elementary tools}

The analogue of the following idea is due to Whitley \cite{whitley}.

Let $X$ be a real Banach space and $K\subseteq X$ be a separable weakly
compact set. Denote $J:X\rightarrow X^{**}$ the usual canonical embedding.

\begin{theorem}
$\overline{\limfunc{co}}\left( K\right) $ is weakly compact (the closure is
taken in norm).
\end{theorem}

\begin{bizonyitas}
By Alaoglu's theorem, the unit ball $\mathcal{B}_{\left( C\left(
K,w^{*}\right) \right) ^{*}}$ of the Banach space $\left( C\left(
K,w^{*}\right) \right) ^{*}$ is $w^{*}$-compact. The restriction operator $%
R:X^{*}\rightarrow C\left( K\right) ,$ $\varphi \mapsto \varphi _{\mid K}$
is continuous, so 
\[
M:=\left\{ \Phi R:\Phi \in \mathcal{B}_{\left( C\left( K\right) \right)
^{*}}\right\} \subseteq X^{**}
\]
is also convex and $w^{*}$-compact, hence $\overline{\limfunc{co}}\left(
J\left( K\right) \right) \subseteq M.$ Now it is enough to prove that every $%
\Phi R\in M$ $\left( \Phi \in \mathcal{B}_{\left( C\left( K\right) \right)
^{*}}\right) $ belonging to the $w^{*}$-closure of $\overline{\limfunc{co}}%
\left( J\left( K\right) \right) ,$ is an element of $J\left( X\right) .$
Suppose it is not; then by Hahn-Banach theorem we obtain a functional $A\in 
\emph{B}_{X^{***}}$ such that $AJ=0$ and $\alpha :=A\left( \Phi R\right) >0.$
Let $\left( x_{n}\right) $ be a norm dense sequence in $K.$ By Goldstine's
theorem we get for every integer $n$ that there is a functional $\varphi
_{n}\in \mathcal{B}_{X^{*}}$ such that 
\[
\left| \Phi R\varphi _{n}-\alpha \right| ,\left| \varphi _{n}x_{1}\right|
,\left| \varphi _{n}x_{2}\right| ,\ldots ,\left| \varphi _{n}x_{n}\right|
\leq \frac{1}{n}\;.
\]
Then the sequence $\left( \varphi _{n}\right) $ tends to $0$ at every $x_{k},
$ thus, $\left( R\varphi _{n}\right) $ being bounded tends to $0$ pointwise
on the whole $K.$ Hence by the Remark \ref{remark} after Theorem \ref
{convexifi} we get that $\Phi R\varphi _{n}\rightarrow 0.$ On the other
hand, by the choice, $\Phi R\varphi _{n}$ obviously tends to $\alpha ,$
which is a contradiction.
\end{bizonyitas}

\end{document}